\documentclass[12pt]{article}
\usepackage{amsmath,amssymb,amsthm,euscript}

\makeatletter
\@addtoreset{equation}{section}

\newcommand{\gs}[1]{\textbf{#1}}
\newcommand{\ds}{\displaystyle}
\newcommand{\fr}[2]{\frac{#1}{#2}}
\newcommand{\cd}{\cdot}

\renewcommand{\l}{\left}
\renewcommand{\r}{\right}

\newcommand{\vsb}{\vspace{2mm}}

\newcommand{\q}{\quad}
\newcommand{\qq}{\qquad}

\newcommand{\maru}[1]{{\ooalign{\hfil#1\/\hfil\crcr
\raise.167ex\hbox{\mathhexbox20D}}}}

\newcommand{\la}{\langle}
\newcommand{\ra}{\rangle}

\newcommand{\abs}[1]{\lvert{#1}\rvert}

\DeclareMathOperator*{\tensor}{\otimes}
\DeclareMathOperator*{\fusion}{\boxtimes}

\newcommand{\Z}{\mathbb{Z}}
\newcommand{\C}{\mathbb{C}}

\newcommand{\N}{\mathbb{N}}
\newcommand{\Q}{\mathbb{Q}}

\newcommand{\aut}{\mathrm{Aut}}

\renewcommand{\hom}{\mathrm{Hom}}

\newcommand{\id}{\mathrm{id}}
\newcommand{\ind}{\mathrm{Ind}}

\newcommand{\pii}{\pi \sqrt{-1}\, }

\newcommand{\vacuum}{\mathrm{1\hspace{-3.2pt}l}}
\newcommand{\vac}{\vacuum}

\newcommand{\pf}{\gs{Proof:}\q}

\renewcommand{\S}{\mathcal{S}}
\newcommand{\indw}{\ind_{V^0}^{V_D}(W,\varphi)}
\newcommand{\red}{\mathrm{Red}}

\theoremstyle{plain}
\newtheorem{thm}{Theorem}[section]
\newtheorem{prop}[thm]{Proposition}
\newtheorem{lem}[thm]{Lemma}
\newtheorem{cor}[thm]{Corollary}

\theoremstyle{definition}
\newtheorem{df}[thm]{Definition}

\theoremstyle{remark}
\newtheorem{rem}[thm]{Remark}

\author{Hiroshi Yamauchi \vsb\\
  {\it \small Graduate School of Mathematics},\\
  {\it \small University of Tsukuba, Ibaraki 305-8571, Japan}\\
  {\small E-mail: {\sf hirocci@math.tsukuba.ac.jp}}
}

\title{Module categories of simple current extensions of \\
vertex operator algebras}

\date{}

\begin{document}

\baselineskip 6mm

\maketitle

\begin{abstract}
  We study module categories of simple current extensions of 
  rational $C_2$-cofinite vertex operator algebras of CFT-type 
  and prove that they are semisimple.
  We also develop a method of induced modules for the simple 
  current extensions.
\end{abstract}


\section{Introduction}

Simple current extensions of vertex operator algebras (VOAs) 
seem to be the most realistic and effective construction of 
new vertex operator algebras, and have many applications.
It is shown in \cite{H4} that the famous 
moonshine VOA $V^\natural$ is a $\Z_2$-graded simple current 
extension of the charge conjugate orbifold $V_\Lambda^+$ of 
the lattice VOA $V_\Lambda$ associated to the Leech lattice 
$\Lambda$.
On the other hand, it is initiated in \cite{DMZ} to study
the moonshine VOA as a module for a sub VOA inside 
$V^\natural$.
The merit of this method is that we can define some 
automorphisms using fusion rules for a sub VOA, see \cite{M1}.
In \cite{LLY} they constructed a series of $\Z_2$-graded
simple current extensions of the unitary Virasoro VOAs 
and computed all the fusion rules for the extensions.
A new kind of symmetry is obtained, and the symmetry 
is shown to be a 3A-triality of the Monster in \cite{SY}.
In \cite{LYY}, they also use simple current extensions 
to study $V^\natural$.
Thus, simple current extensions are useful to study 
the Monster through $V^\natural$. 
However, to define automorphisms by using a sub VOA, 
we have to establish the rationality of the sub VOA.
In this paper, we establish the rationality of the 
simple current extensions of rational $C_2$-cofinite vertex 
operator algebras of CFT-type.

The main tool we use in this paper is the "associativity".
By Huang's results, intertwining operators among simple 
current modules have a nice property so that
in the representation theory of a simple current extension 
with an abelian symmetry we can find certain twisted algebras 
associated to pairs of the abelian group and its orbit spaces.
The twisted algebras can be considered as a deformation or a
generalization of group rings and play a powerful role in our
theory.
Using the twisted algebras, we can show that every module 
for a simple current extension of a rational $C_2$-cofinite 
VOA is completely reducible.
Furthermore, we can parameterize irreducible modules for 
extensions by irreducible representations of the twisted 
algebras.

We also develop a method of induced modules.
We show that every irreducible module of a rational 
$C_2$-cofinite vertex operator algebra of CFT-type can be 
lifted to a twisted module for a simple current extension 
with an abelian symmetry.
This result concerns the following famous conjecture: 
"For a simple rational vertex operator algebra $V$ and its 
finite automorphism group $G$, the $G$-invariants $V^G$, called 
the $G$-orbifold of $V$, is also a simple rational vertex 
operator algebra.
Moreover, every irreducible module for the $G$-orbifold 
$V^G$ is contained in a $g$-twisted $V$-module for 
some $g\in G$."
Our result is the converse of this conjecture in a sense.
Actually, we prove the following.
Let $V^0$ be a simple rational VOA of CFT-type and $D$ a finite 
abelian group. 
Then a $D$-graded simple current extension 
$V_D= \oplus_{\alpha \in D}V^\alpha$ is $\sigma$-regular 
for all $\sigma\in D^*$.
Here the term "$\sigma$-regular" means that every weak 
$\sigma$-twisted module is completely reducible (cf. \cite{Y}).
Moreover, for an irreducible $V^0$-module $W$, we can attach
the group representation $\chi : D\to \Z/n\Z$ such that the 
powers of $z$ in a $V^0$-intertwining operator of type 
$V^\alpha \times (V^\beta \fusion_{V^0} W)\to 
V^{\alpha +\beta} \fusion_{V^0} W$ are contained in 
$\chi (\alpha) +\Z$ for all $\alpha,\beta\in D$.
Finally we prove that $W$ can be lifted to an irreducible
$\hat{\chi}$-twisted $V_D$-module, where $\hat{\chi}$ is an 
element in $D^*$ defined as $\hat{\chi} (\alpha)= 
e^{-2\pii \chi (\alpha)}$ for $\alpha \in D$.

This paper is organized as follows.
In Section 2 we mainly study "decompositions" of $V_D$-modules as 
$V^0$-modules. It is proved that $V_D$ is $\sigma$-regular for
all $\sigma\in D^*$.
In Section 3 we study "induction" of $V_D$-modules from
$V^0$-modules, which is a converse step of Section 2.
As an application, we complete the classification of all 
$\Z_3$-twisted modules for the vertex operator algebra 
constructed in \cite{SY}.

\section{Representation theory of simple current extensions} 

Throughout this paper, $V^0$ denotes a simple rational 
$C_2$-cofinite vertex operator algebra (VOA) of CFT-type.
In particular, $V^0$ is regular by \cite{ABD}.
In this paper, we assume that the module category of $V^0$ 
is well-known in a sense that all irreducible modules and
all fusion rules are classified.

We recall a definition of tensor product.
A theory of tensor products of modules of a vertex operator
algebra was first developed by Huang and Lepowsky in 
\cite{HL1}--\cite{HL4}, and Li considered it in a formal 
calculus approach in \cite{Li1}.
Huang and Lepowsky's approach contains not only intertwining
operators, but also intertwining maps.
However, since we will work over a simple situation, 
in this paper we use Li's definition which is enough 
for us.

\begin{df}(\cite{Li1}) 
  Let $W^1$, $W^2$ be $V^0$-modules.
  A {\it tensor product} (or a {\it fusion product}) for 
  the ordered pair $(W^1,W^2)$ 
  is a pair $(W^1\fusion_{V^0}W^2 ,F(\cd,z))$ consisting 
  of a $V^0$-module $W^1 \fusion_{V^0} W^2$ and a
  $V^0$-intertwining operator $F(\cd,z)$ of type $W^1\times 
  W^2\to W^1\fusion_{V^0} W^2$ satisfying the following 
  universal property:
  For any $V^0$-module $U$ and any intertwining operator 
  $I(\cd,z)$ of type $W^1\times W^2 \to U$, there exists a
  unique $V^0$-homomorphism $\psi$ from $W^1\fusion_{V^0} 
  W^2$ to $U$ such that $I(\cd,z)=\psi F(\cd,z)$.
\end{df}

\begin{rem}
  It is shown in \cite{HL1}--\cite{HL4} and \cite{Li1} that
  if a VOA $V$ is rational, then a tensor product for 
  any two $V$-modules always exists.
  It follows from definition that if a tensor product exists, then 
  it is unique up to isomorphism.
\end{rem}

For simplicity, we will often denote a tensor product 
$(W^1\fusion_{V^0} W^2, F(\cd,z))$ simply by 
$W^1\fusion_{V^0} W^2$.

\begin{df}
  An irreducible $V^0$-module $U$ is called a {\it simple 
  current} if it satisfies: For any irreducible $V^0$-module 
  $W$, a fusion product $U\fusion_{V^0} W$ is also irreducible.
\end{df}

We study an extension of $V^0$ by simple current modules.
Let $D$ be a finite abelian group and assume that we have a set 
of irreducible simple current $V^0$-modules $\{ V^\alpha \mid 
\alpha \in D\}$ indexed by $D$. 

\begin{df}
  An extension $V_D=\oplus_{\alpha \in D} V^\alpha$ of 
  $V^0$ is called a {\it $D$-graded simple current 
  extension\footnote{We assume that the vacuum vector 
  and the Virasoro vector of $V_D$ is the same as those of 
  $V^0$.}} if $V_D$ carries a structure of a simple VOA such 
  that $Y(x^\alpha,z)x^\beta\in V^{\alpha+\beta}((z))$ for 
  any $x^\alpha \in V^\alpha$ and $x^\beta \in V^\beta$. 
\end{df}

\begin{rem}
  It is shown in Proposition 5.3 of \cite{DM2} that the 
  VOA-structure of $V_D$ over $\C$ is unique.
\end{rem}

In the following, $V_D$ always denotes a $D$-graded 
simple current extension of $V^0$.
Let $D^*$ be the dual group of $D$, that is, the group of 
characters of $D$.
With the canonical action, $D^*$ faithfully acts on $V_D$.
So we may view $D^*$ as a subgroup of $\aut (V_D)$.
Then it follows from \cite{DM1} that $V^\alpha$ and $V^\beta$ 
are inequivalent $V^0$-modules if $\alpha \ne \beta$.

\begin{lem}\label{C_2}
  $V_D$ is $C_2$-cofinite.
\end{lem}

\pf
Since we have assumed that $V^0$ is a $C_2$-cofinite VOA of 
CFT-type, we can use a method of spanning sets for modules 
developed by several authors in \cite{B} \cite{M3} \cite{Y},
and we see that all $V^\alpha$, $\alpha\in D$, are 
$C_2$-cofinite $V^0$-modules.
Then the assertion immediately follows.
\qed

\subsection{Twisted algebra $A_\lambda(D,\S_W)$} 

Let $W$ be an irreducible $V^0$-module.
In this subsection we describe a construction of the twisted 
algebra from $W$ and $V^\alpha$, $\alpha\in D$.
Since all $V^\alpha$, $\alpha \in D$, are simple current 
$V^0$-modules, all $V^\alpha \fusion_{V^0} W$, $\alpha \in D$,  
are also irreducible $V^0$-modules.
By the results of Huang \cite{H1} \cite{H2} \cite{H3}, 
fusion products among $V^0$-modules satisfy the associativity.
Therefore $V^\alpha \fusion_{V^0} W \ne 0$ for all $\alpha \in D$ 
and $D_W:=\{ \alpha \in D \mid V^\alpha \fusion_{V^0} 
W \simeq W\}$ forms a subgroup of $D$.
Set $\S_W :=D/D_W$. Then $\S_W$ naturally admits an 
action of $D$.
By definition, $D_W$ acts on $\S_W$ trivially.
Let $s\in \S_W$ and take a representative $\alpha \in D$ 
such that $s=\alpha+D_W$.
We should note that irreducible $V^0$-modules 
$V^\alpha\fusion_{V^0} W$ and $V^\beta \fusion_{V^0} W$ 
are isomorphic if and only if $\alpha -\beta \in D_W$.
Thus, the equivalent class of $V^\alpha \fusion_{V^0}W$ is 
independent of choice of a representative $\alpha$ in 
$\alpha +D_W$ and hence determined uniquely.
So for each $s\in \S_W$, we define $W^s:= V^\alpha \fusion_{V^0} W$ 
after fixing a representative $\alpha \in D$ such that 
$s=\alpha+D_W$.

Let $\alpha,\beta,\gamma \in D$ and $s\in \S_W$.
It follows from the associativity of fusion products that
$V^\alpha \fusion_{V^0} W^s=W^{s+\alpha}$, where
$s+\alpha$ denotes the action of $\alpha \in D$ to 
$s\in \S_W$.
Take basis $I^\alpha_s(\cd,z)$ of the 1-dimensional spaces 
of $V^0$-intertwining operators of type $V^\alpha \times 
W^s \to W^{s+ \alpha}$.
By an associative property of $V^0$-intertwining operators 
(cf. \cite{H1}--\cite{H3}), there are (non-zero) scalars 
$\lambda_s(\alpha,\beta)\in \C$ such that the following 
equality holds:
$$
  \la\nu, I^\alpha_{s+ \beta} (x^\alpha,z_1) 
    I^\beta_s(x^\beta,z_2) w^s \ra
  = \lambda_s(\alpha,\beta) \la \nu, I^{\alpha+\beta}_s
    (Y_{V_D}(x^\alpha,z_0)x^\beta,z_2)w^s\ra |_{z_0=z_1-z_2},
$$
where $x^\alpha \in V^\alpha$, $x^\beta\in V^\beta$, 
$w^s\in W^s$ and $\nu \in (W^{s+\alpha +\beta})^*$.
We normalize intertwining operators $I^0_s(\cd,z)$ to satisfy
$I^0_s(\vac ,z)=\id_{W^s}$.
In other words, $I^0_s(\cd,z)$ are vertex operators on 
$V^0$-modules $W^s$.
By considering 
$$
  \la \nu', I^\alpha_{s+\beta +\gamma}
  (x^\alpha,z_1) I^\beta_{s+ \gamma}(x^\beta,z_2)
  I^\gamma_s(x^\gamma,z_3)w^s\ra ,
$$
we can deduce a relation
$$
  \lambda_{s+ \gamma}(\alpha,\beta)
    \lambda_s(\alpha+\beta,\gamma)
  = \lambda_s(\alpha,\beta+\gamma)
    \lambda_s(\beta,\gamma).
$$
By the normalization $I_s^0(\vac,z)=\id_{W^s}$, 
the $\lambda_s(\cd,\cd)$'s above also satisfy a condition 
$\lambda_s(0,\alpha) =\lambda_s(\alpha,0)=1$ for all 
$\alpha \in D$ and $s\in \S_W$.
Using $\lambda_s(\alpha,\beta)$, we introduce the twisted 
algebra.
Let $q(s)$, $s\in \S_W$, be formal symbols and $\C \S_W 
:= \oplus_{s\in \S_W} \C q(s)$ a linear space spanned by them.
We define a multiplication on $\C \S_W$ by $q(s)\cd q(t):= 
\delta_{s,t}q(s)$. 
Then $\C \S_W$ becomes a semisimple commutative associative 
algebra isomorphic to $\C^{\oplus \abs{\S_W}}$.
Let $U(\C \S_W):= \{ \sum_{s\in \S_W} \mu_s q(s) \mid 
\mu_s\in \C^*\}$ be the set of units in $\C \S_W$.
Then $U(\C \S_W)$ forms a multiplicative group in $\C \S_W$.
Define an action of $\alpha \in D$ on $\C \S_W$ by 
$q(s)^\alpha := q(s-\alpha)$.
Then $U(\C \S_W)$ is a multiplicative right $D$-module.
Set $\bar{\lambda}(\alpha,\beta)=\sum_{s\in \S_W} 
\lambda_s(\alpha,\beta)^{-1} q(s) \in U(\C \S_W)$.
Then $\bar{\lambda}(\cd,\cd)$ defines a 2-cocycle 
$D\times D \to U(\C \S_W)$ because it satisfies a 2-cocycle
condition
$$
  \bar{\lambda}(\alpha,\beta)^\gamma\cd 
  \bar{\lambda}(\alpha+\beta,\gamma) 
  = \bar{\lambda}(\alpha,\beta+\gamma) \cd
  \bar{\lambda}(\beta,\gamma).
$$
Since the space of $V^0$-intertwining operators of 
type $V^\alpha\times W^s \to W^{s+\alpha}$ is 1-dimensional, 
$\bar{\lambda}$ is unique up to 2-coboundaries.
Namely, $\bar{\lambda}$ defines an element of
the second cohomology group $H^2(D,U(\C \S_W))$.
Let $\C [D]=\oplus_{\alpha\in D} \C e^\alpha$ be the group 
ring of $D$ and set
$$
  A_\lambda(D,\S_W):= \C [D]\tensor_\C \C \S_W =\l\{ \sum
  \mu_{\alpha,s} e^\alpha\tensor q(s) \mid \alpha \in D, 
  s \in \S_W, \mu_{\alpha,s} \in \C \r\}
$$
and define a multiplication $*$ on $A_\lambda (D,\S_W)$ by 
$$
  e^\alpha \tensor q(s)* e^\beta \tensor q(t)
  := \lambda_t(\alpha,\beta)^{-1} e^{\alpha+\beta}\tensor 
  q(s)^\beta \cd q(t).
$$
Then $A_\lambda(D,\S_W)$ equipped with the product $*$ forms
an associative algebra with the unit element 
$\sum_{s\in S_W} e^0\tensor q(s)$.
We call $A_\lambda (D,\S_W)$ the {\it twisted algebra} 
associated to a pair $(D,\S_W)$.

\begin{rem}
  The algebra $A_\lambda (D,\S_W)$ is called the 
  {\it generalized twisted double} in \cite{DY}.
  It naturally appears in the orbifold theory, and has been 
  studied in many papers.
  For reference, see \cite{DVVV} \cite{DM3} \cite{DY} 
  \cite{Ma}.
\end{rem}

Take an $s\in \S_W$ and set $\C [D]\tensor q(s):= 
\oplus_{\alpha \in D} \C e^\alpha \tensor q(s)$.
Then $\C [D]\tensor q(s)$ is a subalgebra of 
$A_\lambda (D,\S_W)$.
It has a subalgebra $\C [D_W] \tensor q(s):= 
\oplus_{\alpha \in D_W} \C e^\alpha \tensor q(s)$ which is 
isomorphic to the twisted group algebra $\C^{\lambda_s} [D_W]$ 
of $D_W$ associated to a 2-cocycle $\lambda_s(\cd,\cd)^{-1} 
\in Z^2(D_W,\C)$.
There is a one-to-one correspondence between the category of 
$\C^{\lambda_s} [D_W]$-modules and the category of 
$A_\lambda (D,\S_W)$-modules given as below:

\begin{thm}\label{simple modules}
  (\cite{Ma} \cite[Theorem 3.5]{DY})
  The functors
  $$
  \begin{array}{llll}
    \ind_{\C^{\lambda_s}[D_W]}^{A_\lambda (D,\S_W)}
    : & M \in \C^{\lambda_s}[D_W]\text{-}\mathrm{Mod} 
    & \mapsto & \ds
    \C [D]\tensor q(s) \hspace{-4.5mm} 
    \bigotimes_{\C [D_W]\tensor q(s)} \hspace{-5mm} M
    \in A_\lambda (D,\S_W)\text{-}\mathrm{Mod},
    \vsb\\
    \red_{\C^{\lambda_s}[D_W]}^{A_\lambda (D,\S_W)}
    : & N\in A_\lambda (D,\S_W)\text{-}\mathrm{Mod} & \mapsto &
    e^0\tensor q(s) N\in \C^{\lambda_s}[D_W]
    \text{-}\mathrm{Mod}
  \end{array}
  $$
  define equivalences between the module categories 
  $\C^{\, \lambda_s} [D_W]$-$\mathrm{Mod}$ and $A_\lambda 
  (D,\S_W)$-$\mathrm{Mod}$.
  In particular, $A_\lambda (D,\S_W)$ is a semisimple algebra.
\end{thm}

\subsection{Untwisted Modules} 

Let $M$ be an indecomposable weak $V_D$-module.
Since $V^0$ is regular, we can find an irreducible 
$V^0$-submodule $W$ of $M$.
We use the same notation for $D_W$, $\S_W$, $A_\lambda(D,\S_W)$
and $\C^{\lambda_s}[D_W]$ as previously.
We should note that the definition of $D_W$ is 
independent of the choice of an irreducible component $W$.
One can show the following.

\begin{lem}\label{sy}
  (\cite[Lemma 3.6]{SY})
  All $V^\alpha \cd W=\{ \sum a_n w \mid a\in V^\alpha, w\in W, 
  n\in \Z\}$, $\alpha \in D$, are non-trivial irreducible 
  $V^0$-submodules of $M$.
\end{lem}

Denote by $D=\bigsqcup_{i=1}^n (t^i+D_W)$ a coset decomposition 
of $D$ with respect to $D_W$.
Set $V_{D_W+t^i}:=\oplus_{\alpha \in D_W} V^{\alpha +t^i}$.
Then $V_{D_W}$ forms a sub VOA of $V_D$, which is a $D_W$-graded
simple current extension of $V^0$, 
and $V_D=\oplus_{i=1}^n V_{D_W+t^i}$ forms a $D/D_W$-graded
simple current extension of $V_{D_W}$.
As we have assumed that $M$ is an indecomposable $V_D$-module, 
every irreducible $V^0$-submodule is isomorphic to one of 
$V^{t^i} \fusion_{V^0} W$, $i=1,\dots,n$.

\begin{rem}\label{D/D_W}
  Denote by $M_{D_W+t^i}$ the sum of all irreducible 
  $V^0$-submodules of $M$ isomorphic to $V^{t^i}
  \fusion_{V^0} W$.
  Then we have the following decomposition of $M$ into a 
  direct sum of isotypical $V^0$-components with a 
  $D/D_W$-grading:
  $$
    M=\bigoplus_{t=1}^n M_{D_W+t^i},\q  
    V_{D_W+t^i}\cd M_{D_W+t^j}= M_{D_W+t^i+t^j}.
  $$
  In particular, if $M$ is irreducible under $V_D$, then 
  each $M_{D_W+t^i}$ is an irreducible $V_{D_W}$-module.
  Thus, viewing $V_D$ as a $D/D_W$-graded simple current 
  extension of $V_{D_W}$,  we can regard $M$ as a 
  $D/D_W$-stable $V_D$-module 
  (for the $D/D_W$-stability, see Definition 
  \ref{stability} below).
\end{rem}

For each $s\in \S_W$, we denote $V^s\fusion_{V^0} W$ by 
$W^s$ by abuse of notation (because it is well-defined).
Since all $V^\alpha$, $\alpha\in D$, are simple current 
$V^0$-modules, there are unique $V^0$-intertwining operators 
$I^\alpha_s (\cd,z)$ of type $V^\alpha \times W^s\to 
W^{s+\alpha}$ up to scalar multiples.
We choose $I^0_s (\cd,z)$ to satisfy the condition 
$I^0_s(\vac,z) = \id_{W^s}$, i.e., $I^0_s(\cd,z)$ defines 
the vertex operator on a $V^0$-module $W^s$ for each $s\in 
\S_W$.
Then, by Huang \cite{H1}--\cite{H3}, there exist scalars 
$\lambda_s (\alpha,\beta)\in \C$ such that
\begin{equation}\label{associativity}
  \la \nu, I^\alpha_{s+\beta} (x^\alpha,z_1) 
    I^\beta_s (x^\beta,z_2) w\ra
  = \lambda_s (\alpha,\beta)
  \la \nu, I^{\alpha +\beta}_s(Y_{V_D}
  (x^\alpha,z_0) x^\beta,z_2)w\ra |_{z_0=z_1-z_2},
\end{equation}
where $x^\alpha \in V^\alpha$, $x^\beta \in V^\beta$, 
$w\in W^s$ and $\nu \in (W^{s+\alpha+\beta})^*$.
Then by the same procedure as in the previous subsection, 
we can find the 2-cocycle $\bar{\lambda}(\cd,\cd)\in 
H^2(D,U(\C \S_W))$ and construct the twisted algebra 
$A_\lambda(D,\S_W)$.
By assumption, $M$ is a direct sum of some copies of $W^s$, 
$s\in \S_W$, as a $V^0$-modules so that we have $M \simeq 
\oplus_{s\in \S_W} W^s \tensor \hom_{V^0}(W^s,M)$.
Set $U^s:= \hom_{V^0}(W^s,M)$.
Clearly, all of $U^s$, $s\in \S_W$, are not zero because of 
Lemma \ref{sy}.
On $W^s\tensor U^s$, the vertex operator of 
$x^\alpha \in V^\alpha$ can be written as
\begin{equation}\label{vertex operator}
  Y_M(x^\alpha,z)|_{W^s\tensor U^s}
  =I^\alpha_s(x^\alpha,z)\tensor \phi_s(\alpha)
\end{equation}
with some $\phi_s (\alpha) \in \hom_\C (U^s,U^{s+\alpha})$.
Using \eqref{associativity} we can show that
$$
  \phi_{s+\beta}(\alpha) \phi_s (\beta)=\lambda_s 
  (\alpha,\beta)^{-1} \phi_s(\alpha+\beta)
$$
and hence we can define an action of $A_{\lambda}(D,\S_W)$ on 
$\oplus_{s\in \S_W}U^s$ by 
$e^\alpha\tensor q(s)\cd \mu:= \delta_{s,t} \phi_s(\alpha) \mu$ 
for $\mu \in U^t$. 

\begin{lem}
  Under the action above, $\oplus_{s\in \S_W} U^s$ becomes an 
  $A_\lambda (D,\S_W)$-module.
\end{lem}

Then we have

\begin{prop}\label{structure}
  Suppose that $M$ is irreducible under $V_D$.
  Then $U^s$ is an irreducible $\C^{\lambda_s}[D_W]$-module 
  for every $s\in \S_W$.
  Moreover, $\oplus_{s\in \S_W} U^s$ is an irreducible 
  $A_\lambda (D,\S_W)$-module.
\end{prop}

\pf
Since $W^s\tensor U^s$ is an irreducible $V_{D_W}$-submodule 
of $M$, $U^s$ is an irreducible $\C^{\lambda_s}[D_W]$-module.
Moreover, there is a canonical $A_\lambda(D,\S_W)$-homomorphism
from $P=A_\lambda (D,\S_W)\tensor_{\C [D_W]\tensor q(s)} U^s
=\{ \C [D]\tensor q(s)\} \tensor_{\C [D_W]\tensor q(s)} U^s$
to $\oplus_{s\in \S_W}U^s$. 
As 
$$
  P=\bigoplus_{t\in \S_W} e^t\tensor q(s) 
  \tensor_{\C [D_W]\tensor q(s)} U^s
$$ 
and is irreducible under 
$A_\lambda (D,\S_W)$ by Theorem \ref{simple modules}, 
we see that $U^t= e^t\tensor q(s) \tensor U^s$ 
and hence $P= \oplus_{s\in \S_W} U^s$.
Consequently, $\oplus_{s\in \S_W} U^s$ is also irreducible 
under $A_\lambda (D,\S_W)$.
\qed

\begin{cor}
  If $M$ is an irreducible $V_D$-module, then irreducible 
  $V^0$-components $W^s$ and $W^t$ have the same 
  multiplicity in $M$ for all $s,t\in \S_W$.
\end{cor}

\pf 
Because $U^s$ and $U^t$ have the same dimension.
\qed

\begin{thm}\label{decomposition}
  Let $V^0$ be a rational $C_2$-cofinite vertex operator algebra
  of CFT-type, and let $V_D=\oplus_{\alpha\in D} V^\alpha$ be 
  a $D$-graded simple current extension of $V^0$.
  Then an indecomposable $V_D$-module $M$ is completely 
  reducible under $V_D$. Consequently, $V_D$ is regular.
  As a $V^0$-module, an irreducible $V_D$-submodule of $M$ has
  the shape $\oplus_{s\in \S_W}  W^s \tensor U^s$ with each 
  $U^s$ an irreducible $\C [D_W]\tensor q(s) \simeq 
  \C^{\lambda_s} [D_W]$-module for all $s\in \S_W$.
  Moreover, all $U^t$, $t\in \S_W$, are determined by 
  one of them, say $U^s$, by the following rule: 
  $$
    U^t \simeq \red^{A_\lambda (D,\S_W)}_{\C^{\lambda_t}
      [D_W]}
    \ind_{\C^{\lambda_s}[D_W]}^{A_\lambda (D,\S_W)} U^s.
  $$
\end{thm}

\pf
Since $M=\oplus_{s\in \S_W} W^s\tensor \hom_{V^0}(W^s,M)$ and
the space $\oplus_{s\in \S_W} \hom_{V^0}(W^s,M)$ carries a 
structure of a module for a semisimple algebra 
$A_\lambda(D,\S_W)$ by Proposition \ref{structure}, 
$M$ is also a completely reducible $V_D$-module because of 
\eqref{vertex operator}.
Since $V^0$ is regular, all weak $V_D$-modules are completely 
reducible.
So $V_D$ is also regular.
Now assume that $M$ is an irreducible $V_D$-module.
The decomposition is already shown.
It remains to show that $U^t$ is determined by $U^s$ by 
the rule as stated.
It is shown in the proof of Proposition \ref{structure} 
that $U^t =e^t\tensor q(s) \tensor U^s$.
It is easy to see that $e^t\tensor q(s) \tensor U^s =
\red^{A_\lambda (D,\S_W)}_{\C^{\lambda_t}[D_W]}
\ind_{\C^{\lambda_s}[D_W]}^{A_\lambda (D,\S_W)} U^s$.
The proof is completed.
\qed

By the theorem above and Theorem \ref{simple modules}, 
the number of inequivalent irreducible $V_D$-modules containing
$W$ as a $V^0$-submodule is $\dim_\C Z(\C^{\lambda_s} [D_W])$.
In particular, if $D_W=0$, then the structure of a $V_D$-module
containing $W$ is uniquely determined by its $V^0$-module 
structure.
This fact is already observed in Proposition 3.8 of \cite{SY}.
For convenience, we introduce the following notion.

\begin{df}\label{stability}
  An irreducible $V_D$-module $N$ is said to be 
  {\it $D$-stable} 
  if $D_W=0$ for some irreducible $V^0$-submodule $W$.
\end{df}

It is obvious that the definition of the $D$-stability is 
independent of the choice of an irreducible $V^0$-submodule $W$.
Let $N^i$, $i=1,2,3$, be irreducible $D$-stable $V_D$-modules,
and let $W^i$ be irreducible $V^0$-submodules of $N^i$ for 
$i=1,2,3$.
Set $W^{i,\alpha}:= V^\alpha \fusion_{V^0} W^i$.
Then $W^{i,\alpha}\simeq W^{i,\beta}$ as $V^0$-modules 
if and only if $\alpha =\beta$, 
and $N^i$ as a $V^0$-module is isomorphic to 
$\oplus_{\alpha \in D} W^{i,\alpha}$.
The following lemma is a simple modification of Lemma 3.12 of 
\cite{SY}.
We expect that this lemma would be used in our future study.

\begin{lem}\label{fusion}
  Let $N^i=\oplus_{\alpha\in D} W^{i,\alpha}$ be as above.
  Then we have the following isomorphism:
  $$
    \binom{N^3}{N^1\ N^2}_{V_D} \simeq 
    \bigoplus_{\alpha\in D}\binom{W^{3,\alpha}}{W^1\ W^2}_{V^0}.
  $$
  In the above, $\binom{N^3}{N^1\ N^2}_{V_D}$ denotes the
  space of $V_D$-intertwining operators of type $N^1\times N^2
  \to N^3$.
\end{lem}

\pf 
By the Huang's result in \cite{H3}, we can replace the 
assumption on the Virasoro element in Lemma 3.12 of \cite{SY}
by the $C_2$-cofinite condition.
\qed

\begin{rem}
  Let $M$ be an irreducible $V_D$-module and $W$  an 
  irreducible $V^0$-submodule of $M$.
  Even if $D_W\ne 0$, we can apply the lemma above to $M$.
  We may consider $V_D$ as a $D/D_W$-graded simple 
  current extension of $V_{D_W}$ as in Remark \ref{D/D_W}.
  Then we can view $M$ as a $D/D_W$-stable $V_D$-module.
  So by replacing $D$ by $D/D_W$, we can apply the lemma above
  to $M$.
\end{rem}

\subsection{Twisted modules} 

Let $\sigma$ be an automorphism on $V_D$ such that $V^0$ is 
contained in $V_D^{\la \sigma\ra}$, the space of 
$\sigma$-invariants of $V_D$.
Denote by $V_D^{(r)}$ subspaces $\{ a\in V_D \mid \sigma a 
=e^{2\pii r/\abs{\sigma}} a \}$ for $0\leq r\leq 
\abs{\sigma}-1$, where $\abs{\sigma}$ is order of $\sigma$.
Then $V_D^{(r)}$ are $V^0$-submodules of $V_D$ so that there is 
a partition $D=\bigsqcup_{i=0}^{\abs{\sigma}-1} D^{(i)}$ 
such that $V_D^{(i)}=\oplus_{\alpha \in D^{(i)}} V^\alpha$.
One can easily verify that $D^{(0)}$ is a subgroup of $D$ and 
each $D^{(i)}$ is a coset of $D$ with respect to $D^{(0)}$.
Namely, if $V^0\subset V_D^{\la \sigma\ra}$, then $\sigma$ is 
identified with an element of $D^*$, the dual group of $D$.
Conversely, it is clear from definition that $D^*$ is a 
subgroup of $\aut (V_D)$.
Thus, we have

\begin{lem}
  An automorphism $\sigma\in \aut (V)$ satisfies $V^0\subset 
  V_D^{\la \sigma \ra}$ if and only if $\sigma \in D^*$.
\end{lem}

The lemma above tells us that an automorphism $\sigma$ is 
consistent with the $D$-grading of $V_D$ if and only if $\sigma$ 
belongs to $D^*$.
We consider $\sigma$-twisted $V_D$-modules.
Let $M$ be an indecomposable admissible $\sigma$-twisted 
$V_D$-module.
By definition, there is a decomposition
$$
  M=\bigoplus_{i=0}^{\abs{\sigma}-1} M^{(i)}
$$
such that $V_D^{(i)}\cd M^{(j)}\subset M^{(i+j)}$.
It is obvious that each $M^{(i)}$ is a $V^0$-module.
Let $W$ be an irreducible $V^0$-submodule of $M^{(0)}$, and
let $D_W$, $\S_W$, $A_\lambda(D,\S_W)$ and $\C^{\lambda_s}[D_W]$
be as in Section 2.1.
By replacing $M$ by $\sum_{\alpha \in D} V^\alpha \cd W$ if
necessary, we may assume that all $V^\alpha\fusion_{V^0} W$, 
$\alpha \in D_W$, are contained in $M^{(0)}$ so that $D_W$ is a 
subgroup of $D^{(0)}$.
Since $M$ is a completely reducible $V^0$-module, we have the 
following decomposition:
$$
  M= \bigoplus_{s\in \S_W} W^s\tensor \hom_{V^0}(W^s,M),
$$
where we set $W^s:= V^s\fusion_{V^0} W$ for $s\in \S_W$ by 
abuse of notation.
Set $U^s:= \hom_{V^0}(W^s,M)$ for $s\in \S_W$.
As we did before, we can find a 2-cocycle 
$\bar{\lambda} \in H^2(D,U(\C \S_W))$ and a representation 
of $A_\lambda (D,\S_W)$ on the space 
$\oplus_{s\in \S_W} U^s$.
Thus, by the same argument, we can show the following.

\begin{thm}
  Let $\sigma \in D^*(\subset \aut (V_D))$.
  Viewing as a $V^0$-module, an indecomposable admissible 
  $\sigma$-twisted $V_D$-module $M$ has the shape
  $$
    M=\bigoplus_{s\in \S_W} W^s\tensor U^s
  $$
  such that the space $\oplus_{s\in \S_W} U^s$ carries a 
  structure of an $A_\lambda (D,\S_W)$-module.
  In particular, $M$ is a completely reducible $V_D$-module.
  If $M$ is irreducible under $V_D$, then each $U^s$, 
  $s\in \S_W$, is irreducible under $\C [D_W]\tensor q(s)$, and
  also $\oplus_{s\in \S_W} U^s$ is irreducible under 
  $A_\lambda(D,\S_W)$.
  Moreover, for each pair $s$ and $t\in \S_W$, $U^s$ and $U^t$ is 
  determined by the following rule:
  $$
    U^t \simeq \red_{\C^{\lambda_t}[D_W]}^{A_\lambda (D,\S_W)}
    \ind_{\C^{\lambda_s}[D_W]}^{A_\lambda (D,\S_W)} U^s.
  $$
  Hence, all $W^s$, $s\in \S_W$, have the same multiplicity 
  in $M$.
\end{thm}

\begin{rem}
  Since $D_W\subset D^{(0)}$, we note that the decomposition 
  above is a refinement of the decomposition $M=
  \oplus_{i\in \Z/\abs{\sigma}\Z} M^{(i)}$.
\end{rem}

By the theorem above, $V_D$ is $\sigma$-rational for all $\sigma
\in D^*$. 
More precisely, we can prove that $V_D$ is $\sigma$-regular, 
that is, every weak $\sigma$-twisted $V_D$-module is 
completely reducible (cf.{} \cite{Y}).

\begin{cor}
  An extension $V_D$ is $\sigma$-regular for all $\sigma\in D^*$.
\end{cor}

\pf
Let $M$ be a weak $\sigma$-twisted $V_D$-module.
Take an irreducible $V^0$-submodule $W$ of $M$, which is possible
because $V^0$ is regular.
Then $\sum_{\alpha \in D}V^\alpha\cd W$ is a $\sigma$-twisted 
admissible $V_D$-submodule.
As we have shown that $V_D$ is $\sigma$-rational, 
$\sum_{\alpha \in D} V^\alpha \cd W$ is a completely reducible 
$V_D$-module.
Thus, $M$ is a sum of irreducible $V_D$-submodules and hence 
$M$ is a direct sum of irreducible $V_D$-submodules.
\qed

\section{Module categories of simple current extensions} 

In the previous section, we studied a representation theory of a 
$D$-graded simple current extension $V_D$.
In this section, we proceed a theory of induced modules, which is 
given as a converse step of the previous section.
As a result, we can complete the classification of irreducible 
modules for $V_D$.

\subsection{Induced modules} 

Let $W$ be an irreducible $V^0$-module.
We define the stabilizer $D_W$, the orbit space 
$\S_W$, intertwining operators $I^\alpha_s(\cd,z)$, where 
$\alpha \in D$ and $s\in \S_W$, the twisted algebra 
$A_\lambda(D,\S_W)$ and the twisted group ring 
$\C^{\lambda_s}[D_W]$ as in Section 2.1.
We set $W^s:= V^s\fusion_{V^0} W$ for $s\in \S_W$, as we did
previously.
Denote by $h(s)$ the top weight of a $V^0$-module $W^s$, which
is a rational number by Theorem 11.3 of \cite{DLiM3}.
It follows from definition that the powers of $z$ in 
an intertwining operator $I^\alpha_s(\cd,z)$ are contained 
in $h(\alpha +s)-h(s) +\Z$.
We set $\chi(\alpha ,s):= h(\alpha +s)-h(s)\in \Q$.
The following assertion is crucial for us.

\begin{lem}\label{key}
  The following hold for any $\alpha,\beta\in D$ and $s\in 
  \S_W$:\\
  {} \q 
  (i) $\chi (\alpha, \beta+s) -\chi (\alpha,s)\in \Z\,;$ 
  \q
  (ii) $\chi (\alpha,\beta +s)+\chi(\beta,s)-\chi 
  (\alpha+\beta,s)\in \Z$.
\end{lem}

\pf
Recall that the following results are established by Huang 
in \cite{H1}--\cite{H3}:
\begin{align}
  & 
  \la \nu, I^\alpha_{s+\beta}(x^\alpha,z_1) 
    I^\beta_s(x^\beta,z_2) w^s\ra
  = \epsilon_s(\alpha,\beta) 
    \la \nu, I^\beta_{s+\alpha}(x^\beta,z_2) 
    I^\alpha_s(x^\alpha,z_1) w^s\ra, \label{eq1}
  \vsb\\
  & 
  \la \nu, I^\alpha_{s+\beta}(x^\alpha,z_1) 
    I^\beta_s(x^\beta,z_2) w^s\ra
  = \lambda_s(\alpha,\beta) 
    \la \nu, I^{\alpha+\beta}_s(Y_{V_D}(x^\alpha,z_0) 
    x^\beta,z_2) w^s\ra |_{z_0=z_1-z_2}, \label{eq2}
\end{align}
where $x^\alpha \in V^\alpha$, $x^\beta\in V^\beta$, 
$w^s \in W^s$, $\nu\in (W^{s+\alpha+\beta})^*$, 
$\epsilon_s(\alpha,\beta)$ is a suitable scalar in $\C$, 
and the equals above mean that the left hand side and 
the right hand side are analytic extensions 
of each other.
Since all $I^\alpha_s(\cd,z)$ are intertwining operators 
among modules involving simple currents, we note that by 
the convergence property in \cite{H1}--\cite{H3} 
the right hand side of \eqref{eq2} has the form 
$
  z_2^{\chi (\alpha+\beta,s)+r} z_0^s f_1(z_0/z_2) 
  |_{z_0=z_1-z_2}
$ 
in the domain $\abs{z_2}>\abs{z_1-z_2}>0$, 
where $r$ and $s$ are some integers and $f_1(x)$ is an
analytic function on $\abs{x}<1$.
Therefore, we have
\begin{equation}\label{eq3}
\begin{array}{l} 
  (z_1-z_2)^N \la \mu,I^\alpha_{s+\beta}(x^\alpha,z_1)
  I^\beta_s(x^\beta,z_2) w^2\ra
  \vsb\\
  = \epsilon_s(\alpha,\beta) (z_1-z_2)^N
  \la \mu, I^\alpha_{s+\beta} (x^\beta,z_2)
  I^\beta_s(x^\alpha,z_2) w^2\ra
\end{array}
\end{equation}
in the domain $\abs{z_1}>\abs{z_1-z_2}\geq 0$, 
$\abs{z_2}>\abs{z_1-z_2}\geq 0$ for sufficiently large 
$N$.
Since $z^{-\chi(\alpha,s)} I^\alpha_s(x^\alpha,z) w^s$ 
contains only integral powers of $z$, both 
$$
  z_1^{-\chi (\alpha,\beta +s)} z_2^{-\chi (\beta,s)} 
  \la \nu, I^\alpha_{s+\beta}(x^\alpha,z_1) 
  I^\beta_s(x^\beta,z_2) w^s\ra
$$ 
and 
$$
  z_1^{-\chi (\alpha,s)} z_2^{-\chi (\beta,s+\alpha)} 
  \la \nu, I^\alpha_{s+\alpha} (x^\beta,z_2) 
  I^\beta_{s+\alpha}(x^\beta,z_2) w^s\ra
$$
contain only integral powers of $z_1$ and $z_2$.
Thus, by \eqref{eq1}, \eqref{eq3} and the convergence property 
in \cite{H1}--\cite{H3}, 
we obtain the following equality of the meromorphic functions:
$$
\begin{array}{l}
  (z_1-z_2)^N \, 
    \iota_{12}^{-1} z_1^{-\chi(\alpha,s+\beta)} 
    z_2^{-\chi (\beta,s)} \la \nu, I^\alpha_{s+\beta} 
    (x^\alpha,z_1) I^\beta_s(x^\beta,z_2) w^s\ra 
  \vsb\\
  = (z_1-z_2)^N \epsilon_s(\alpha,\beta)
    \iota_{21}^{-1} z_1^{-\chi (\alpha ,s)} 
    z_2^{-\chi(\beta,s+\alpha)} 
    \cd z_1^{\chi(\alpha,s)-\chi(\alpha,s+\beta)} 
    z_2^{\chi(\beta,s+\alpha)-\chi(\beta,s)}
  \vsb\\
  \qq \times  
    \la \nu,I^\beta_{s+\alpha}(x^\beta,z_2) I^\alpha_s
    (x^\alpha,z_1) w^s\ra .
\end{array}
$$
Since the equality above holds for any choices of $\log z_1$ and 
$\log z_2$ in the definitions of $z_1^r=e^{r\log z_1}$ and 
$z_2^r=e^{r\log z_2}$ (cf. \cite{H1}--\cite{H3}), 
we have $\chi(\alpha,s)-\chi(\alpha,s+\beta)\in \Z$ 
and $\chi(\beta,s+\alpha)-\chi(\beta,s)\in \Z$.
This proves (i).
The proof of (ii) is similar.
By \eqref{eq2} and the convergence property in 
\cite{H1}--\cite{H3}, 
we obtain the following equality of the meromorphic functions:
$$
\begin{array}{l}
  \lambda_s(\alpha,\beta)^{-1} 
    \iota_{12}^{-1} 
    z_1^{-\chi (\alpha,s+\beta)} z_2^{-\chi(\beta,s)}
    \la \nu, I^\alpha_{s+\beta}(x^\alpha,z_1)
    I^\beta_s(x^\beta,z_2) w^s\ra
  \vsb\\
  = \iota^{-1}_{20} z_2^{-\chi(\alpha+\beta,s)}
    \la \nu, I^{\alpha+\beta}_s(Y_{V_D}(x^\alpha,z_0)
    x^\beta,z_2) w^s \ra 
    \cd (z_2+z_0)^{-\chi(\alpha,s+\beta)}
    z_2^{\chi(\alpha+\beta,s)- \chi(\beta,s)}
    |_{z_0=z_1-z_2}.
\end{array}
$$
Again the equality holds for any choices of $\log (z_1-z_2)$ and 
$\log z_2$ in the definitions of $(z_1-z_2)^r=e^{r\log (z_1-z_2)}$ 
and $z_2^r=e^{r\log z_2}$ (cf. \cite{H1}--\cite{H3}).
Since $(z_2+z_0)^r=z_2^r(1+z_0/z_2)^r$ and $(1+x)^r=\sum_{i\geq 0}
\binom{r}{i} x^i$ is analytic in the domain $\abs{x} <1$, 
we see that $\chi(\alpha+\beta,s)-\chi(\alpha,s+\beta)-
\chi (\beta,s) \in \Z$. 
This completes the proof of (ii).
\qed

By the lemma above, we find that $\chi(\alpha,s)+\Z$ is 
independent of $s\in \S_W$.
So we may set $\chi(\alpha):= \chi (\alpha,s)$ for 
$\alpha \in D$.
Then by (ii) of Lemma \ref{key} we find that $\chi(\cd)$ 
satisfies the homomorphism condition $\chi(\alpha+\beta)+\Z=
\chi(\alpha) +\chi(\beta)+\Z$.
Since $\chi(\alpha)\in \Q$, there exists an $n\in \N$ such that 
$\chi$ defines a group homomorphism from $D/D_W$ to
$\Z/n\Z=\{ j/n +\Z \mid 0\leq j\leq n-1\}$.
It is clear that $\chi$ naturally defines an element 
$\hat{\chi}$ of $D^*$ by $\hat{\chi}(\alpha) := 
e^{-2\pii \chi (\alpha)}$.
Thus $\chi$ gives rise to an element of $\aut (V_D)$.
In the following, we will construct irreducible 
$\hat{\chi}$-twisted $V_D$-modules which contain $W$ as
$V^0$-submodules.

Take an $s\in S_W$.
Let $\varphi$ be an irreducible representation of
$\C^{\lambda_s}[D_W]$ on a space $U$. 
For each $t\in S_W$, set 
$$
  U^t:= \red_{\C^{\lambda_t} [D_W]}^{A_\lambda(D,S_W)} 
  \ind_{\C^{\lambda_s} [D_W]}^{A_\lambda(D,S_W)} U.
$$
Then each $U^t$, $t\in \S_W$, is a 
$\C^{\lambda_t} [D_W]$-module and a direct sum 
$\oplus_{s\in S_W} U^s$ naturally (and uniquely) 
carries a structure of an irreducible $A_\lambda(D,S_W)$-module.
Set 
$$
  \ind_{V^0}^{V_D} (W,\varphi)
  := \bigoplus_{s\in S_W} W^s\tensor U^s
$$
and define the vertex operator $\hat{Y}(\cd,z)$ of $V_D$ on 
$\ind_{V^0}^{V_D} (W,\varphi)$ by
$$
  \hat{Y}(x^\alpha,z) w^t\tensor \mu^t
  := I^\alpha_t(x^\alpha,z)w^t\tensor \{ e^\alpha\tensor 
  q(t)\cd \mu^t\}
$$
for $x^\alpha\in V^\alpha$ and $w^t\tensor \mu^t \in 
W^t\tensor U^t$.
We prove

\begin{thm}\label{induced module}
  $(\ind_{V^0}^{V_D}(W,\varphi), \hat{Y}(\cd,z))$ is an 
  irreducible $\hat{\chi}$-twisted $V_D$-module.
\end{thm}

\pf
Since the powers of $z$ in $\hat{Y}(x^\alpha,z)$ are contained 
in $\chi(\alpha)+\Z$, we only need to show the commutativity 
and the $\hat{\chi}$-twisted associativity of vertex operators.
We use a technique of generalized rational functions developed
in \cite{DL}.
Let $x^\alpha \in V^\alpha$, $x^\beta \in V^\beta$, $w^s\tensor 
\mu^s \in W^s\tensor U^s$ and $\nu\in \indw^*$.
We note that $z^{-\chi (\alpha)} I^\alpha_s(x^\alpha,z) w^s\in 
W^{s+\alpha}((z))$.
For sufficiently large $N\in \N$, we have
$$
\begin{array}{l}
  (z_1-z_2)^N
    \iota_{12}^{-1}  z_1^{-\chi (\alpha)} z_2^{-\chi (\beta)}
    \la \nu, \hat{Y}(x^\alpha,z_1) 
    \hat{Y}(x^\beta,z_2) w^s\tensor \mu^s\ra
  \vsb\\
  = (z_1-z_2)^N \iota_{12}^{-1}
    z_1^{-\chi(\alpha)} z_2^{-\chi (\beta)} 
    \la \nu, I^\alpha_{s+\beta}(x^\alpha,z_1) 
    I^\beta_s(x^\beta,z_2) w^s 
  \vsb\\
  \qq\q \tensor \{ e^\alpha\tensor q(s+\beta) \cd 
    (e^\beta\tensor q(s) \cd \mu^s)\}\ra
  \vsb\\
  = (z_1-z_2)^N \iota_{123}^{-1}
    z_1^{-\chi(\alpha)} z_2^{-\chi (\beta)}
    \la \nu, I^\alpha_{s+\beta}(x^\alpha,z_1) 
    I^\beta_s(x^\beta,z_2)
    I^0_s(\vac,z_3)w^s
  \vsb\\
  \qq\q  \tensor \{ \lambda_s(\alpha,\beta)^{-1} 
    e^{\alpha +\beta}\tensor q(s)\cd \mu^s\} \ra
  \vsb\\
  = \iota_{345}^{-1} (z_3+z_4)^{-\chi(\alpha)} 
    (z_3+z_5)^{-\chi(\beta)}
    \vsb\\
  \qq\q 
    \cd \la \nu, \lambda_s(\alpha,\beta)
    I^{\alpha+\beta}_s((z_4-z_5)^N 
    Y_{V_D}(x^\alpha,z_4) Y_{V_D}(x^\beta,z_5) 
    \vac,z_3)w^s
    \vsb\\
  \qq\q 
    \tensor \{ \lambda_s(\alpha,\beta)^{-1} 
    e^{\alpha+\beta} \tensor q(s) \cd \mu^s\}\ra 
    |_{z_4=z_1-z_3,\, z_5=z_2-z_3}
\end{array}
$$
$$
\begin{array}{l}
  = \iota_{354}^{-1} (z_3+z_4)^{-\chi (\alpha)} 
    (z_3+z_5)^{-\chi (\beta)} 
    \la \nu, I^{\alpha+\beta}_s((z_4-z_5)^N 
    Y_{V_D}(x^\beta,z_5) Y_{V_D}(x^\alpha,z_4) 
    \vac,z_3)w^s
  \vsb\\
  \qq\q \tensor \{ e^{\alpha+\beta} 
    \tensor q(s) \cd \mu^s\}\ra |_{z_4=z_1-z_3,\, z_5=z_2-z_3}
  \vsb\\
  = (z_1-z_2)^N \iota_{123}^{-1} z_1^{-\chi(\alpha)} 
    z_2^{-\chi(\beta)} 
    \la \nu, \lambda_s(\beta,\alpha)^{-1}
    I^\beta_{s+\alpha} (x^\beta,z_2) I^\alpha_s(x^\alpha,z_1)
    I^0_s(\vac,z_3)w^s
  \vsb\\
  \qq\q \tensor \{ (\lambda_s(\beta,\alpha)
    e^\beta\tensor q(s+\alpha)*e^\alpha\tensor 
    q(s))\cd \mu^s\} \ra
  \vsb\\
  = (z_1-z_2)^N \iota_{12}^{-1} z_1^{-\chi(\alpha)} 
    z_2^{-\chi(\beta)} 
    \la \nu, \hat{Y}(x^\beta,z_2) 
    \hat{Y}(x^\alpha,z_1) w^s\tensor \mu^s\ra.
\end{array}
$$ 
Therefore, we get the commutativity.
Similarly, we have 
$$
\begin{array}{l}
  \iota_{12}^{-1}
    z_1^{-\chi (\alpha)+N} z_2^{-\chi (\beta)}
    \la \nu, \hat{Y}(x^\alpha,z_1) \hat{Y}(x^\beta,z_2) 
    w^s\tensor \mu^s \ra
  \vsb\\
  = \iota_{12}^{-1} 
    z_1^{-\chi (\alpha)+N} z_2^{-\chi (\beta)}
    \la \nu, I^\alpha_{s+\beta} (x^\alpha,z_1)
    I^\beta_s(x^\beta,z_2) w^s \tensor 
    \{ e^\alpha\tensor q(s+\beta)\cd (e^\beta \tensor q(s)\cd \mu^s)\} 
    \ra
  \vsb\\
  = \iota_{20}^{-1} \la \nu, \lambda_s(\alpha,\beta)
    (z_2+z_0)^{-\chi (\alpha)+N} z_2^{-\chi (\beta)}
    I^{\alpha +\beta}_s(Y_{V_D}(x^\alpha,z_0)x^\beta,z_2) w^s
  \vsb\\
    \qq\q \tensor \{ \lambda_s(\alpha,\beta)^{-1} 
    e^{\alpha+\beta} \tensor q(s)\cd \mu^s\} 
    \ra |_{z_0=z_1-z_2}
  \vsb\\
  = \iota_{20}^{-1} \la \nu, (z_2+z_0)^{-\chi (\alpha)+N}
    z_2^{-\chi (\beta)} \hat{Y}(Y_{V_D}(x^\alpha,z_0)x^\beta,z_2)
    w^s\tensor \mu^s\ra |_{z_0=z_1-z_2}.
\end{array}
$$
Hence, we obtain the associativity.
\qed

Suppose that a simple VOA $V$ and a finite group $G$ acting on 
$V$ is given.
Then the $G$-invariants $V^G$ of $V$, called the $G$-orbifold 
of $V$, is also a simple VOA by \cite{DM1}.
It is an important problem to classify the module category of 
$V^G$ in the orbifold conformal field theory.
It was conjectured in \cite{DVVV} that every irreducible 
$V^G$-module appears in a $g$-twisted $V$-module for some 
$g\in G$.
In our case, $V^0$ is exactly the $D^*$-invariants of the 
extension $V_D$.
By Theorem \ref{induced module}, 
we see that the conjecture is true for a pair $(V_D,D^*)$.

\begin{thm}
  Let $V^0$ be a rational, $C_2$-cofinite and CFT-type
  VOA, and $D$ a finite abelian group.
  Let $V_D=\oplus_{\alpha \in D} V^\alpha$ be a
  $D$-graded simple current extension of $V^0$.
  Then every irreducible $V^0$-module $W$ is contained in 
  an irreducible $\sigma$-twisted $V_D$-module for some 
  $\sigma\in D^*$.
  Moreover, $\sigma$ is uniquely determined by $W$.
\end{thm}

\subsection{Application to the 3A-algebra for the Monster}

Here we give an application of Theorem \ref{induced module}.
We use the notation of \cite{SY} without any comments.
A typical example of a vertex operator algebra generated by
two conformal vectors is studied in \cite{SY}.
Denote it by $U$ as in \cite{SY}.
It has a shape $U=U^0\oplus U^1\oplus U^2$ with 
$U^0= W(0) \tensor N(0)$, 
$U^1= W(\fr{2}{3})^+ \tensor N(\fr{4}{3})^+$, 
$U^2= W(\fr{2}{3})^- \tensor N(\fr{4}{3})^-$, 
and it is shown in \cite{SY} that $U$ is a $\Z_3$-graded 
simple current extension of $U^0$.
Define $\zeta\in \aut (U)$ by $\zeta|_{U^i}=e^{2\pii i/3}\cd 
\id_{U^i}$.
Then by Theorem \ref{induced module} every irreducible 
$U^0$-module uniquely is induced to be a 
$\zeta^j$-twisted $U$-module as follows:
\vsb\\
\begin{tabular}{lll}
  (i) untwisted $U$-modules:
  & $\ind_{U^0}^U W(h)\tensor N(k)$, 
  & $h\in \{ 0,\fr{2}{5}\}$, 
    $k\in \{ 0,\fr{1}{7}, \fr{5}{7}\}$.
\vsb\\
  (ii) $\zeta$-twisted $U$-modules:
  & $\ind_{U^0}^U W(h)\tensor N(k)^+$, 
  & $h\in \{ 0,\fr{2}{5}\}$, 
    $k\in \{ \fr{4}{3}, \fr{10}{21}, \fr{1}{21}\}$.
\vsb\\
  (iii) $\zeta^2$-twisted $U$-modules:
  & $\ind_{U^0}^U W(h)\tensor N(k)^-$, 
  & $h\in \{ 0,\fr{2}{5}\}$, 
    $k\in \{ \fr{4}{3},\fr{10}{21},\fr{1}{21}\}$.
\end{tabular}
\vsb\\
In the above, $\ind_{U^0}^U X$ denotes 
$X \oplus (U^1\fusion_{U^0} X) \oplus (U^2\fusion_{U^0}X)$ for 
an irreducible $U^0$-module $X$ (note that $U^i\fusion_{U^0} X
\ne U^j\fusion_{U^0}X$ if $i\not\equiv j \mod 3$).
It is shown in \cite{SY} that $\aut (U)=S_3$.
Thus, we have a classification of $\Z_3$-twisted $U$-modules.

\begin{thm}
  Every irreducible $U^0$-module is uniquely lifted to an
  irreducible $\zeta^i$-twisted $U$-module for some $i=0,1,2$.
\end{thm}

\begin{rem}
  It is shown in \cite{SY} that $U$ is contained in the famous
  moonshine VOA $V^\natural$ \cite{FLM}, and that $\zeta$ 
  naturally induces a 3A element of Monster, the automorphism
  group of $V^\natural$.
\end{rem}


\small
\baselineskip 12pt

\begin{thebibliography}{ABiCD} 


\bibitem[ABD]{ABD}
  T. Abe, G. Buhl and C. Dong, Rationality, regularity, and 
  $C_2$-cofiniteness, preprint, math.QA/0204021.

\bibitem[B]{B}
  G. Buhl, A spanning set for VOA modules, {\it J. Algebra} 
  {\bf 254} (2002), 125--151.

\bibitem[DVVV]{DVVV}
  R. Dijkgraaf, C. Vafa, E. Verlinde and H. Verlinde, 
  The operator algebra of orbifold models, {\it Comm. Math. 
  Phys.} {\bf 123} (1989), 485--526.

\bibitem[DL]{DL}
  C. Dong and J. Lepowsky, Generalized vertex algebras and 
  relative vertex operators, Progress in Math. {\bf 112},
  Birkh\"{a}user, Boston, 1993. 

\bibitem[DLiM1]{DLiM1}
  C. Dong, H. Li and G. Mason, Hom functor and the associativity 
  of tensor products of modules for Vertex operator algebras,
  {\it J. Algebra} {\bf 188} (1997), 443--475.

\bibitem[DLiM2]{DLiM2}
  C. Dong, H. Li and G. Mason, Simple currents and extensions of
  vertex operator algebras, {\it Comm. Math. Phys.} {\bf 180}, 
  671--707.

\bibitem[DLiM3]{DLiM3}
  C. Dong, H. Li and G. Mason, Modular-invariance of trace
  functions in orbifold theory and generalized Moonshine, 
  {\it Comm. Math. Phys.} {\bf 214} (200), 1--56.

\bibitem[DM1]{DM1}
  C. Dong and G. Mason, On quantum Galois theory, 
  {\it Duke Math. J.} {\bf 86} (1997), 305--321.

\bibitem[DM2]{DM2}
  C. Dong and G. Mason, Rational vertex operator algebras and 
  the effective central charge, preprint, math.QA/0201318.

\bibitem[DM3]{DM3}
  C. Dong and G. Mason, Vertex operator algebras and Moonshine;
  A Survey, {\it Advanced Studies in Pure Mathematics} {\bf 24}, 
  Progress in Algebraic Combinatorics, Mathematical Society of 
  Japan.

\bibitem[DMZ]{DMZ}
  C. Dong, G. Mason and Y. Zhu, Discrete series of the Virasoro
  algebra and the moonshine module, Proc. Symp. Pure. Math.,
  American Math. Soc. {\bf 56} II (1994), 295--316.


\bibitem[DY]{DY}
  C. Dong and G. Yamskulna, Vertex operator algebras, generalized
  doubles and dual pairs, {\it Math. Z.} {\bf 241} (2002), 
  397--423.


\bibitem[FLM]{FLM}
  I.B. Frenkel, J. Lepowsky and A. Meurman, Vertex Operator 
  Algebras and the Monster, Academic Press, New York, 1988.



\bibitem[H1]{H1}
  Y.-Z. Huang, A theory of tensor products for module categories 
  for a vertex operator algebra, IV {\it J. Pure Appl. Algebra} 
  {\bf 100} (1995), 173--216.

\bibitem[H2]{H2}
  Y.-Z. Huang, Virasoro vertex operator algebras, (non-meromorphic)
  operator product expansion and the tensor product theory, 
  {\it J. Algebra} {\bf 182} (1996), 201--234.

\bibitem[H3]{H3}
  Y.-Z. Huang, Differential equations and intertwining operators,
  math.QA/0206206.

\bibitem[H4]{H4}
  Y.-Z. Huang, A nonmeromorphic extension of the moonshine
  module vertex operator algebra, {\it in} "Moonshine, the 
  Monster and Related Topics, Proc. Joint Summer Research
  Conference, Mount Holyoke, 1994" 
  (C. Dong and G. Mason, Eds.), 
  Contemporary Math., pp. 123--148, Amer. Math. Soc., 
  Providence, RI, 1996.

\bibitem[HL1]{HL1}
  Y.-Z. Huang and J. Lepowsky, 
  Toward a theory of tensor product for representations of 
  a vertex operator algebra, in Proc. 20th Intl. Conference
  on Diff. Geom. Methods in Theoretical Physics, New York,
  1991, ed. S. Catto and A. Rocha, World Scientific, 
  Singapore, 1992, Vol. 1, 344-354.

\bibitem[HL2]{HL2}
  Y.-Z. Huang and J. Lepowsky, 
  Tensor products of modules for a vertex operator algebra 
  and vertex tensor categories, in: Lie Theory and Geometry, 
  ed. R. Brylinski, J.-L. Brylinski, V. Guillemin, V. Kac,
  Birkhauser, Boston, 1994, 349--383.

\bibitem[HL3]{HL3}
  Y.-Z. Huang and J. Lepowsky, 
  A theory of tensor products for module categories for a 
  vertex operator algebra, I, {\it Sel. Math.} {\bf 1}
  (1995), 699-756.
\bibitem[HL4]{HL4} 
  Y.-Z. Huang and J. Lepowsky, 
  A theory of tensor products for module categories for a 
  vertex operator algebra, III, {\it J. Pure Appl. Alg.}
  {\bf 100} (1995), 141-171.



\bibitem[La]{L}
  C.H. Lam, Induced modules for orbifold vertex operator 
  algebras, {\it J. Math. Soc. Japan} {\bf 53} (2001), 541--557.

\bibitem[LLY]{LLY}
  C.H. Lam, N. Lam and H. Yamauchi, Extension of Virasoro vertex
  operator algebra by a simple module, {\it IMRN} {\bf 11} 
  (2003), 577--611.

\bibitem[LYY]{LYY}
  C.H. Lam, H. Yamada and H. Yamauchi, in preparation.

\bibitem[Li]{Li1}
  H. Li, An analogue of the Hom functor and a generalized 
  nuclear democracy theorem, {\it Duke Math. J.} {\bf 93} 
  (1998), 73--114.


\bibitem[Ma]{Ma}
  G. Mason, The quantum double of a finite group and its role 
  in conformal field theory, Proc. LMS Lecture Notes {\bf 212}, 
  CUP.

\bibitem[M1]{M1}
  M. Miyamoto, Griess algebras and conformal vectors in vertex
  operator algebras, {\it J. Algebra} {\bf 179} (1996), 528--548.

 
\bibitem[M2]{M3}
  M. Miyamoto, Modular invariance of vertex operator algebras 
  satisfying \\
  $C_2$-cofiniteness, to appear in {\it Duke. Math. J.}

\bibitem[SY]{SY}
  S. Sakuma and H. Yamauchi, Vertex operator algebra with two
  Miyamoto involutions generating $S_3$, 
  {\it J. Algebra} {\bf 267} (2003), 272--297.


\bibitem[Y]{Y}
  H. Yamauchi, Modularity on vertex operator algebras arising 
  from semisimple primary vectors, math.QA/0209026.


\end{thebibliography}
\end{document}